\documentclass[11pt,a4paper]{article}
\usepackage{times,amssymb,epsfig,rotating}

\newenvironment{keywords}{ \noindent {\small\bf Key Words}:}{ }

\newcommand{\be}{\begin{equation}}
\newcommand{\ee}{\end{equation}}

\newtheorem{theorem}{Theorem}
\newtheorem{lemma}{Lemma}
\newtheorem{corollary}{Corollary}

\newcommand{\bdm}{\begin{displaymath}}
\newcommand{\edm}{\end{displaymath}}

\newcommand{\ba}{\begin{array}}
\newcommand{\ea}{\end{array}}

\begin{document}

\title{Index Information Algorithm with Local Tuning
       for Solving Multidimensional Global Optimization Problems with
       Multiextremal Constraints\thanks{This research was  supported by the
       following grants: FIRB RBNE01WBBB, FIRB RBAU01JYPN, and RFBR
       01--01--00587. }
      }

\author{Yaroslav~D.~Sergeyev \\ DEIS,
         Universit\`{a}
         della Calabria, 87030 Rende (Cs), Italy, and\\
         N.I. Lobachevsky State University of Nizhni Novgorod ,
                  Russia,\\ \it{yaro@si.deis.unical.it} \\
         Paolo~Pugliese \\ DEIS, Universit\`{a} della Calabria,
         87030 Rende (Cs), Italy,\\ \it{pugliese@unical.it} \\
         Domenico~Famularo \\ ICAR--CNR, c/o DEIS, Universit\`{a}
         della Calabria, 87030 Rende (Cs), Italy,\\
         \it{famularo@isi.cs.cnr.it}
         }


\maketitle

\begin{abstract}
Multidimensional optimization problems where the objective
function and the constraints are multiextremal non-differentiable
Lipschitz functions (with unknown Lipschitz constants) and the
feasible  region is a finite collection of robust nonconvex
subregions are considered. Both the objective function and the
constraints may be partially defined. To solve such problems an
algorithm is proposed, that uses Peano space-filling curves and
the index scheme to reduce the original problem to a H\"{o}lder
one-dimensional one. Local tuning on the behaviour of the
objective function and constraints  is used during the work of
the global optimization procedure in order to accelerate the
search. The method neither uses penalty coefficients nor
additional variables. Convergence conditions are established.
Numerical experiments confirm the good performance of the
technique.
\end{abstract}


\section{Introduction}
In the last decades there has been a growing interest in
approaching Global Optimization  problems by different numerical
techniques (see, for example, \cite{Bomze}--\cite{HansenE},
\cite{HorstPard}--\cite{Mockus}, \cite{Pinter,Ratschek},
\cite{SS2000}--\cite{Zhigljavsky} and references given therein).
Such an interest is motivated by a large number of real-life
applications where such problems arise (see, for example,
\cite{BaBo92}--\cite{Blondel},
\cite{Crivelli,FaPuSe99,FloPar,Grossman,HorstPard,Nemirovskii,Pard2000,PardShaXue,%
Pinter,Ritter,Sen,Se99,Shao}). These problems often lead to  deal
with multiextremal non-differentiable objective function and
constraints. In such a context the Lipschitz condition becomes the
unique information about the problem.

It has been proved  by Stephens and Baritompa~\cite{SteBar} that
if the only information about the objective function $\varphi(y)$
is that it belongs to the class of Lipschitz functions and the
Lipschitz constant is unknown, there does not exist any
deterministic or stochastic algorithm that, after a finite number
of evaluations of $\varphi(y)$, is able to provide an
underestimate of the global minimum $\varphi(y^*)$. Of course,
this result is very discouraging because usually in practice it is
difficult to know the constant. Nevertheless, the necessity to
solve practical problems remains. That is why in such problems
instead of the statement (S1) `{\it Find an algorithm able to stop
in a given time and provide an $\varepsilon$-approximation of the
global optimum $\varphi(y^*)$}' the statement (S2) `{\it Find an
algorithm able to stop in  a given time and return the lowest
value of $\varphi(y)$ obtained by the algorithm}' is used. Under
this statement, either the computed solution (possibly improved by
a local analysis) is accepted by final users (engineers,
physicists, chemists, etc.) or the global search is repeated with
changed parameters of the algorithm.

Theoretical analysis of algorithms (depending on parameters) for
solving problems (S2) is similar to analysis of penalty methods.
It is proved for global optimization algorithms that for a {\it
fixed} problem there exists a parameter $P^*$ such that
parameters $P \ge P^*$ allow to solve the problem (obviously, the
parameter $P^*$ is problem-dependent). The problem `{\it How to
determine $P^*$?}' is not discussed in such an analysis, and in
every concrete case it is solved  using additional information
about the problem. For example, in methods using a given
Lipschitz constant (see survey \cite{HanJau}) $P^*$ is  the
Lipschitz constant and it is not discussed how to obtain  it.
Other examples are diagonal methods (see \cite{Pinter}) and
information algorithms (see \cite{Se95,St78}) using similar
parameters. For all those methods it is possible to prove that if
a value $P \ge P^*$ is used as the parameter, then they converge
only to global minimizers. An alternative approach is represented
by methods converging to every point in the search domain (see
\cite{Jones,Torn}).

In this paper  a constrained Lipschitz global minimization problem
is considered. In an informal way it can be stated as follows:

\begin{enumerate}
\item[(i)]
The objective function is multiextremal, non-differentiable,
'black box', and requires a high time to be evaluated;
\item[(ii)]
Constraints are non-convex (or even multiextremal) and
non-differentiable, leading to a complex feasible region
consisting of  disjoint, non-convex sub-regions;
\item[(iii)]
Both the objective function and constraints are Lipschitz
functions with unknown Lipschitz constants;
\item[(iv)]
Both the objective function and constraints may be partially
defined, i.e., if a constraint  is not satisfied at a point, the
rest of constraints and the objective function may be not defined
at that point.
\end{enumerate}

It can be seen from this statement that the problem belongs to the
class of problems considered  in~\cite{SteBar}, thus statement
(S2) will be used hereinafter. One  promising way to face problem
(i)--(iv) is the information approach introduced by Strongin in
\cite{St78,St85,St92}. It uses Peano type space-filling curves
(see~\cite{Bertsimas,Butz,PlatzBar,St78,SS2000}  for examples of
usage of  space-filling curves in mathematical programming) to
reduce the original Lipschitz multi-dimensional problem to a
H\"{o}lder univariate one (a comprehensive presentation  of this
approach can be found in \cite{SS2000}). Global optimization of
H\"{o}lder functions (see \cite{Gourdin,Lera,Vanderbei}) has
given new tools for solving the reduced one-dimensional problem.
Peano curves avoid constructions of support (or auxiliary)
functions usually used in the multi-dimensional Lipschitz
optimization (see, for example,
\cite{HorstPard,Kearfott,Kelley,Pinter,Ratschek,SS2000} and
references given therein). Of course, if the user knows that the
objective function is differentiable and/or the problem is
convex, there is no sense to work with Peano curves and specific
methods explicitly using information about differentiability or
convexity can be applied. In contrast, when you deal with
non-differentiable 'black box' multiextremal problems, it is not
possible to work with sophisticated techniques using derivatives
(or other strong a priori suppositions) and such a reduction to
one dimension can help significantly.

In this paper,  a novel  algorithm belonging to the family of
information methods is proposed. It uses two powerful ideas for
solving problem (i)--(iv). The first one is the index scheme (see
\cite{SeMa95,St78,St85,StronMark,SS2000}), allowing to solve
Lipschitz  problems where both the objective function $\varphi(y)$
and constraints $G_i(y),\, 1\le i\le m,$ may be multiextremal and
partially defined. Its importance increases in this case because
it is not clear how to solve such problems by using, for example,
the penalty approach. In fact, the latter requires that
$\varphi(y)$ and $G_i(y),\, 1\le i\le m,$ are defined over the
whole search domain. It seems that missing values can be simply
filled in with either a big number or the function value at the
nearest feasible point. Unfortunately, in the context of Lipschitz
algorithms, incorporating such ideas can lead to infinitely high
Lipschitz constants, causing degeneration of the methods and
non-applicability of the penalty approach. Thus, for problems
where the Lipschitz condition is almost a unique additional
information, the ability of the index scheme to work with
partially defined problems becomes crucial. Moreover, the index
scheme  does not introduce any additional parameter and/or
variable, whereas the penalty approach requiring determination of
the penalty coefficient.

The second idea being at the basis of the new method is the local
tuning on the behaviour of the objective function  and constraints
(see~\cite{Se95,Se98a,Se99,SeMa95}). Original information methods
work with {\it global} adaptive estimates of the Lipschitz
constants, i.e., the same estimates are used over the {\it whole}
search region for $\varphi(y)$ and functions $G_i(y),\, 1\le i\le
m$. However, global estimates (adaptive or given a priori) of the
Lipschitz constants may provide a  poor information about the
behavior of the objective function over every small sub-region of
$D$. It has been shown in~\cite{Pinter,Se95,Se98a,Se99,SeMa95},
for different classes of global optimization problems, that local
estimates for different sub-regions of $D$ can accelerate the
search significantly. Of course, it is necessary to  balance
global and local information about $\varphi(y)$ obtained by the
method during the search. Such a balancing is very important
because using only the local information can lead to missing the
global solution (see~\cite{Se98b,SteBar}).

The next section contains a formal statement of the problem
(i)--(iv) and presents the new  algorithm. Convergence conditions
of the new method are established in Section~3. Numerical
experiments collected in Section~4  show a satisfactory
performance of the new algorithm in comparison with two global
optimization techniques. Finally, Section~5 gives concluding
remarks.



\section{Theoretical background and the index information\\
algorithm with local tuning}

We start by formulating the Lipschitz optimization problem
satisfying requirements (i)--(iii). Find the constrained global
minimum $\varphi^*$ and at least one  minimizer $y^*$ such that
\begin{equation}
\varphi^* = \varphi(y^*) = \min\{\varphi(y): y\in S, \quad
G_i(y)\le 0,\quad 1\le i\le m\}, \label{f1}
\end{equation}
where the search domain is the hyperinterval
\[
S = \{y\in \mathbb{R}^N: a_j\le y_j\le b_j,\quad 1\le j\le N\},
\]
$\mathbb{R}^N$ is the $N$-dimensional Euclidean space,  and
the          objective function $\varphi(y)$ (henceforth denoted
as $G_{m+1}(y)$) and the functions $G_i(y),\, 1\le i\le m$, of
the constraints can be multiextremal and non-differentiable,
satisfying the Lipschitz condition with unknown constants
$0<\widehat{L}_i<\infty,\,1\le i\le m+1$, i.e.,

\begin{equation}
|G_i(y')-G_i(y'')|\le \widehat{L}_i \, ||y'-y''||,\quad  1\le i\le
 m+1,\quad y',y''\in S. \label{f2}
\end{equation}

Without loss of generality,  we shall consider the search domain
$S=D$, where
\begin{equation}
D=\{y\in \mathbb{R}^N: -2^{-1}\le y_j\le 2^{-1},\quad 1\le j\le
N\}. \label{f3}
\end{equation}
Formulation (\ref{f1})--(\ref{f3}) assumes that all the functions
$G_i(y),\, 1\le i\le m+1,$ can be evaluated in the whole region
$D$. In order to incorporate requirement (iv), we shall assume
that each function $G_i(y)$ is defined and computable only in the
corresponding subset $Q_i\subset D$, where
\begin{equation}
Q_1=D,\quad Q_{i+1}=\{y\in Q_i:G_i(y)\le 0\},\quad 1\le i\le m.
\label{f4}
\end{equation}
The above assumption also imposes the order in which the functions
$G_i(y), 1\le i\le m$, are evaluated. In many applications this
order is determined by the nature of the problem. In other cases
the user introduces a specific order suit for some reasons (for
example,  first verify  easier computable constraints).
 In view of (\ref{f4}), the initial problem
(\ref{f1})--(\ref{f3}) is rewritten as
\begin{equation}
\varphi(y^*)=\min\{G_{m+1}(y): y\in Q_{m+1}\}, \label{f5}
\end{equation}
\begin{equation}
|G_i(y')-G_i(y'')|\le L_i \, ||y'-y''||, \quad y',y''\in Q_i,
\quad 1\le i\le m+1,    \label{f6}
\end{equation}
where $L_i\le\widehat{L}_i$.

\begin{figure}[t]
\framebox[\columnwidth]{\centerline{\psfig{figure=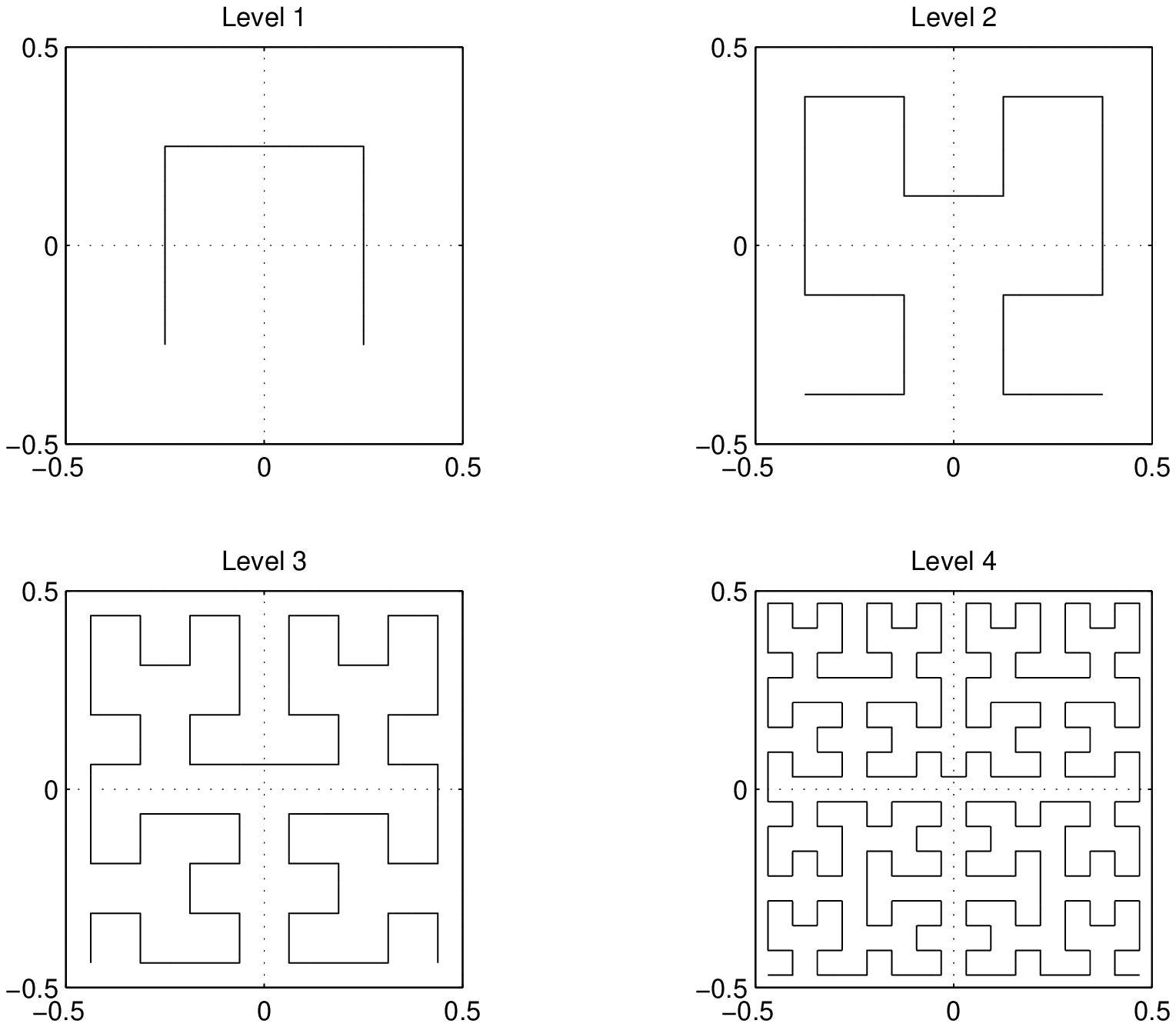,width=0.80\columnwidth,silent=yes}}}
\caption{Approximations of levels one to four to the Peano
(Hilbert) curve in two dimensions.}\label{Fig_1}
\end{figure}

In order to start the description of the method it is necessary to
remind the idea of the space-filling curves. Such curves were
first introduced by Peano in~\cite{Peano} and Hilbert
in~\cite{Hilbert} and are fractal objects constructed on the
principle of self-similarity. They possess the property `to fill'
any cube $D$ in $\mathbb{R}^N$, i.e., they pass through every
point of $D$. An example of construction of such a
 curve (see ~\cite{Sa94,SS2000} for details) is given in Fig.~\ref{Fig_1}.
Naturally, in numerical algorithms, approximations of the curve
are used (Fig.~\ref{Fig_1} presents approximations of levels one
to four).

It has been shown in~\cite{Butz,St78,St92} (see
also~\cite{SS2000}) that the multi-dimensional problem
\begin{equation}
 \varphi^{*} = \varphi(
y^{*} ) = \min \{ \varphi(y):\; y \in D \} , \label {f2.1}
\end{equation}
where $\varphi(y)$ is a Lipschitz function with   constant $L$,
$0 < L < \infty,$ can be reduced by a Peano-type space-filling
curve $y(x)$  to the one-dimensional problem
\begin{equation}
f^*=f(x^*) = \min \{ f(x),\; x \in [ 0,1 ]\} ,  \label {f2.2}
\end{equation}
where the notation $f(x) = \varphi(y(x))$ is used for the obtained
reduced one-dimensional function. Moreover, (see Theorem~8.1 in
\cite{SS2000}), the function $f(x)$ satisfies the H\"{o}lder
condition
\begin{equation}
 | f(x') - f(x'') | \leq H \, | x' - x''|^{1/N}, \,
\quad x',x'' \in [ 0,1], \label {f2.3}
\end{equation}
in the H\"{o}lder metric
\begin{equation}
\rho(x',x'') = | x' - x''|^{1/N}, \label {f2.4}
\end{equation}
where $N$ is from (\ref{f3}) and $H=2\,L\,\sqrt{N+3}$.

By applying the same curve $y(x)$ to the functions $G_i(y),\; 1\le
i\le m+1$, and using designations $g_i(x)=G_i(y(x))$, the
multi-dimensional problem (\ref{f5})--(\ref{f6}) is reduced to
the following constrained one-dimensional problem
(\ref{f2.5})--(\ref{f2.7}):
\begin{equation}
g_{m+1}^*=g_{m+1}(x^*)=\min\{g_{m+1}(x): x\in q_{m+1}\},
                                       \label{f2.5}
\end{equation}
where the region $q_{m+1}$  is defined by the relations
\begin{equation}
q_1=[0,1],\quad q_{i+1}=\{x\in q_i:g_i(x)\le 0\},\quad 1\le i\le
m.                                     \label{f2.6}
\end{equation}
and the reduced functions $g_i(x), 1\le i\le m+1$, satisfy the
corresponding H\"{o}lder conditions
\begin{equation}
|g_i(x')-g_i(x'')|\le H_i \,|x'-x''|^{1/N}, \quad x',x''\in q_i,
\quad 1\le i\le m+1, \label{f2.7}
\end{equation}
with the constants  $H_i=2\,L_i\,\sqrt{N+3}$, where $L_i, 1\le
i\le m+1$, are from (\ref{f6}). Note that for the regions $Q_i$
from (\ref{f4}) and $q_i$ from (\ref{f2.6}) the following relation
holds:
\[
Q_i=\{y(x):x \in q_i\}, \quad 1\le i\le m+1.
\]

This problem may  be rewritten by using the  index scheme
proposed  in \cite{St85,StronMark} (see also \cite{St92}). The
scheme is an alternative to traditional penalty methods. Instead
of combining the objective and constraint functions into a penalty
one (and the need to define in a proper way penalty coefficients),
the scheme does not introduce any additional parameter or
variable and evaluates constraints one at a time at every point
where it has been decided to try to evaluate $g_{m+1}(x)$. The
function $g_i(x)$ is calculated only if all inequalities
\[
g_{j}(x)\le 0,\quad  1\le j<i,
\]
have been satisfied. In its turn, the objective function
$g_{m+1}(x)$ is evaluated only for those points where all the
constraints have been satisfied.

The index scheme juxtaposes to every point of the interval $[0,1]$
an  index
\[
\nu =\nu (x),\quad 1\le \nu \le m+1,
\]
which is defined by the conditions
\begin{equation}
g_{j}(x)\le 0,\quad 1\le j\le \nu -1,\quad g_{\nu }(x)>0,
\label{f2.8}
\end{equation}
where for $\nu =m+1$ the last inequality is omitted. Thus,
$\nu(x)$ represents the number of the first constraint not
satisfied at~$x$ when $\nu(x)<m+1$. Then  $\nu(x) =m+1$ means
that all constraints were satisfied at~$x$ and the objective
function $g_{m+1}(x)$ may be evaluated at this point.

Let us now define an auxiliary function $\Phi (x)$ as follows
\begin{equation}
\Phi (x) = g_{\nu }(x) - \left\{ \begin{array}{ll}
                             0,           & \nu (x)<m+1,\\[5pt]
                             g^{*}_{m+1}, & \nu (x)=m+1,
                                    \end{array}
                            \right.                   \label{f2.9}
\end{equation}
where $g^{*}_{m+1}$ is the solution of the problem
(\ref{f2.5})--(\ref{f2.7}). Due to (\ref{f2.8}) and (\ref{f2.9}),
the function $\Phi (x)$ has the following properties:

\begin{itemize}
\item[i.]
$\Phi (x)>0$, when $\nu (x)<m+1$;
\item[ii.]
$\Phi (x)\ge 0$, when $\nu (x)=m+1$;
\item[iii.]
$\Phi (x)=0$, when $\nu (x)=m+1$ and $g_{m+1}(x)=g^{*}_{m+1}$;
\item[iv.]
$\Phi (x)$ is not continuous at the a priori unknown boundary
points of the sets $q_i,2\le i\le m+1$.
\end{itemize}
Thus, the global minimizer of the constrained problem
(\ref{f2.5})--(\ref{f2.7}) coincides with the so\-lu\-tion
$x^{*}$ of the following unconstrained  problem
\begin{equation}
\Phi (x^{*})=\min \{ \Phi (x): x\in [0,1] \}. \label{f2.10}
\end{equation}
Of course, instead of the fractal Peano curve $y(x)$ its $d$-level
approximation $y_{d}(x)$   is used for evaluation of $\Phi (x)$
(approximations of levels one to four are shown in
Fig.~\ref{Fig_1}, from where it can be seen how a $d$-level
approximation is obtained).

The new index information global optimization algorithm with local
tuning presented below generalizes and evolves two methods. On the
one hand, the information global optimization algorithm with local
tuning proposed in \cite{Se95} for solving the problem
(\ref{f2.1}) and, consequently, the problem
(\ref{f2.2})--(\ref{f2.3}).  On the other hand, the index
algorithm with local tuning proposed in~\cite{SeMa95} for solving
the  one-dimensional problem
\[
\min  \{f(x): x \in [a, b],\; g_i(x) \le 0, \;  1 \le i \le m \}
\]
where both $f(x)$ and  $g_i(x),  1 \le i \le m,$ are Lipschitz
continuous functions.

During every iteration of the algorithm a point $x \in [0,1]$ is
chosen and the value $\Phi (x)$ is evaluated (hereinafter such an
evaluation  will be called a {\it trial} and the corresponding
point $x$ a {\it trial point}).  Suppose now that $k$ iterations
of the algorithm have already been executed (two initial trials
are done at the end points $x^{0} = 0$ and $x^{1} = 1$). The
point $x^{k+1}, k \ge 1$, is determined by the following
algorithm.

\begin{description}
\item[\bf  Step 1.]
The points $x^{0}, \ldots, x^{k}$ of the previous iterations are
renumbered by subscripts as follows
\[ 0 = x_{0} < x_{1} < \cdots < x_i < \cdots  < x_{k} = 1. \]
\item[{\bf Step 2.}]
To each point $x_i$ associate the index $\nu_i=\nu(x_i)$, and the
value
\[
z_i=g_{\nu _i}(x_i) -\left\{ \begin{array}{ll}
                             0,          &  \nu_i<m+1,\\[5pt]
                             z_k^{*},    & \nu_i=m+1,

                                    \end{array}
                            \right.
\]
where
\begin{equation}
z_k^{*}=\min \{ g_{m+1}(x_i): 0 \le i \le k, \nu(x_i)=m+1 \}.
                                     \label{f2.11}
\end{equation}
The value  $z_k^{*}$ estimates the unknown value  $g^{*}_{m+1}$
from~(\ref{f2.9}) on the basis of the available data.

\item[\bf Step 3.]
Calculate lower bounds $\mu_{j}$ for the global H\"older constants
$H_{j}$  of the functions $g_j(x),  1 \le j \le
          m+1, $ as follows
\begin{equation}
  \mu_{j}= \max  \Big\{ \frac{| z_{p}-z_{q} |}{ (x_{p}-x_{q})^{1/N}}:
          0 \le q < p \le k, \nu _{p} = \nu_{q} = j  \Big\}.  \label{f2.12}
\end{equation}
Whenever $\mu_j$ can not be calculated, set $\mu_{j} = 0$.

\item[\bf Step 4.]
For each interval $[x_{i-1},x_i], i = 1, \ldots, k$, calculate the
values
\begin{equation}
M_i= \max \{ \lambda _i, \gamma _i, \xi \} \label{f2.13}
\end{equation}
that estimate the local H\"older constant over the interval
$[x_{i-1},x_i]$. The values $\lambda_i$ and $\gamma_i$ reflect
the influence on $M_i$ of the local and global information
obtained during the previous $k$ iterations, $\xi $ is a small
number - parameter of the method. The values $\lambda_i$ and
$\gamma_i$ are defined by
\begin{equation}
\lambda _i= \max \{ l_i, c_i, r_i\},     \label{f2.14}
\end{equation}
where
\[
\ba{ll} c_i &=  \left\{\!\!\! \begin{array}{cl}
                \displaystyle{\frac{| z_i-z_{i-1} |}{ (x_i-x_{i-1})^{1/N}}}, & \mbox{if} \; \nu_i = \nu_{i-1}   \\[15pt]
                0,    & $ otherwise$
               \end{array}
       \right.\\[25pt]
l_i &= \left\{\!\!\! \begin{array}{cl}
                \displaystyle{\frac{| z_{i-1}-z_{i-2}|}{ (x_{i-1}-x_{i-2})^{1/N}}}, & \mbox{if} \; i\ge 2,\quad  \nu_{i-2} = \nu_{i-1},  \nu_{i-1}\ge\nu_i  \\[15pt]
                0,   &  \mbox{otherwise}
               \end{array}
       \right.\\[25pt]
r_i &= \left\{\!\!\! \begin{array}{cl}
                \displaystyle{\frac{| z_{i+1}-z_i |}{ (x_{i+1}-x_i)^{1/N} }}, & \mbox{if} \; i\le k-1,\quad  \nu_{i+1} = \nu_i,   \nu_i \ge \nu_{i-1}  \\[15pt]
                0,    & $ otherwise$
               \end{array}
       \right.
\ea
\]
and
\begin{equation}
\gamma_i= \frac{\mu_{j}(x_i-x_{i-1})^{1/N} }{X^{\rm max}_{j}},
\hspace{6mm} j=\max\{\nu _i,\nu_{i-1}\},   \label{f2.15}
\end{equation}
where $X^{\rm max}_{j}=  \!\max \{ (x_i-x_{i-1})^{1/N}:%
        \; \max\{\nu_i,\nu_{i-1} \} = j, \quad 1 \le i \le k \}$.

\item[\bf Step 5.]
For each interval $[x_{i-1},x_i], i = 1, \ldots, k$, calculate the
{\it characteristic} of the interval
\begin{equation}
R_i = \left\{ \begin{array}{ll}
                 \Delta_i + \,\displaystyle{ \frac{(z_i-z_{i-1})^{2} }{ (r M_i)^2 \Delta_i }
                    - \, \frac{2\, (z_i+z_{i-1})}{r M_i}}, & \nu_i=\nu_{i-1} \\[10pt]
              2\, \Delta_i - \, \displaystyle{\frac{4\, z_i}{r M_i}}, & \nu_i>\nu_{i-1}\\[10pt]
              2\, \Delta_i - \, \displaystyle{\frac{4\, z_{i-1}}{r M_i}}, & \nu_{i-1}>\nu_i
             \end{array}
       \right.
                             \label{f2.16}
\end{equation}

where $\Delta_i=(x_i-x_{i-1})^{1/N}$ and $r>1$ is a real value --
the reliability parameter of the method.

\item[\bf Step 6.] Choose the interval $[x_{t-1},x_t]$ having the
maximal characteristic as follows
\begin{equation}
t = \min \{ \arg\max \{ R_i: 1 \le i \le k\}
\}.                       \label{f2.17}
\end{equation}

\item[\bf Step 7.] If the interval $[x_{t-1},x_t]$ is such that
\[
(x_t - x_{t-1})^{1/N} > \delta,
\]
where $\delta$ is a given tolerance, go to {\bf Step 8}. Otherwise
{\bf Stop}.

\item[\bf Step 8.]
Execute the $(k+1)$-th iteration at the point
\begin{equation}
x^{k+1} \!= \frac{x_t + x_{t-1}}{2} - \frac{{\rm sign} (z_t -
z_{t-1})|z_t - z_{t-1}|^N}{2\, r\,M_t^N},
                                              \label{f2.18}
\end{equation}
if $\nu_t = \nu_{t-1}$. In all other cases do it at the point
\begin{equation}
 x^{k+1} \!= (x_t + x_{t-1})/2.                       \label{f2.19}
\end{equation}
Set $k=k+1$ and go to {\bf Step 1}.
\end{description}

If trial points with the index $m+1$ have been generated by the
algorithm then the value~$z_k^{*}$ from (\ref{f2.11}) can be taken
as an estimate of the global minimum $\varphi^*$  from (\ref{f1})
and the corresponding point $y_{d}(x_k^{*})$ such that
$z_k^{*}=g_{m+1}(x_k^{*})$ as an estimate of the point $y^*$. If
no points with the index $m+1$ have been generated by the
algorithm then it is necessary to continue the search with
changed parameters of the method.

Let us give a few remarks on the algorithm introduced above. The
information algorithms are derived as optimal statistical decision
functions within the framework of a stochastic model representing
the function to be optimized as a sample of a random function. The
characteristic $R_i$ in terms of the information approach (see
\cite{St78,SS2000}) may be interpreted (after normalization) as
the probability of finding a global minimizer within the interval
$[x_{i-1},x_{i}]$ based on the data available during the current
iteration. The method uses in its work four parameters:  $d, r,
\xi ,$ and $\delta$. Their choice will be discussed in Section~4
while presenting the numerical experiments.

For every sub-interval $[x_{i-1},x_{i}], 1\le i\le k$, {\it
global} estimates $\mu _{j}$ of the {\it global} H\"older
constants $H_{j}, 1\le j\le m+1$, from (\ref{f2.7}) are not used.
In contrast,  {\it local} estimates $M_{i}, 1\le i\le k$, from
(\ref{f2.13}) are adaptively determined. The values $\lambda _{i}$
and $\gamma _{i}$ reflect the influence on $M_{i}$ of the local
and global information obtained during the previous $k$
iterations. When the interval $[x_{i-1},x_{i}]$ is small, then
$\gamma _{i}$ is small too (see (\ref{f2.15})) and, due to
(\ref{f2.13}), the local information represented by $\lambda
_{i}$ has  major importance. The value $\lambda _{i}$ is
calculated by considering the intervals $[x_{i-2},x_{i-1}], \;
[x_{i-1},x_{i}],$ and $ [x_{i},x_{i+1}]$ (see (\ref{f2.14})) as
those which have the strongest influence on the local estimate.
When the interval $[x_{i-1},x_{i}]$ is very wide, the local
information is not reliable and the global information
represented by $\gamma _{i}$ is used. Thus, local and global
information are balanced in the values $M_{i}, 1\le i\le k$. Note
that the method uses the local information over the {\it whole}
search region $[0,1]$ (and, consequently, over the whole
multi-dimensional domain $D$) {\it during} the global search both
for the objective function and constraints (being present in a
implicit form in the auxiliary function $\Phi (x)$).


\section{Convergence conditions}

In the further theoretical consideration it is assumed that each
region $Q_i, 2\le i\le m+1,$ from (\ref{f4}) is a union of $T_i$
disjoint sub-regions $Q_i^j, \; 1\le j\le T_i,$ having positive
volume (in general, the numbers $T_i$ are unknown to the user)
where each sub-region $Q_i^j, 1\le j\le T_i, 2\le i\le m+1,$ is
robust  and can be non-convex. It is assumed that the feasible
region $Q_{m+1}$ is nonempty. Recall that a set $B$ is robust if
for each point $\beta$ belonging to the boundary of $B$ and for
any $ \varepsilon > 0$ there exists an $\varepsilon$-neighborhood
$\varepsilon(\beta)$ such that $\varepsilon(\beta) \bigcap B$ has
a positive volume. Supposition about robustness of the sets
$Q_i^j$ is quite natural in engineering applications (for example,
in optimal control) where the optimal solution must have an
admissible neighborhood of positive volume. This requirement is a
consequence of the inevitable inaccuracy of the physical systems,
and ensures that small changes of the parameters of the system
will not lead the solution to leave the feasible region.

Theoretical results presented in this section have the spirit of
existence theorems. They show (similarly to the penalty approach)
that for a given problem there exist parameters of the method
allowing to solve this {\it fixed} problem. Recall that (see
\cite{SteBar}) the information available from the statement
(\ref{f5})--(\ref{f6}) is not sufficient for establishing concrete
values of the parameters. These should be chosen from an
additional information about the nature of the practical problem
under consideration (see the next section).

The first result links solutions of the original
multi-dimensional problem (\ref{f5})--(\ref{f6}) to the solutions
of the corresponding one-dimensional problem
(\ref{f2.5})--(\ref{f2.7}) reduced via the curve $y_{d}(x)$. It
shows that for any given  $\varepsilon>0$ there  exists a level
$d$ such that the approximation $y_{d}(x)$ passes through the
$\varepsilon$-neighborhood of the global solution $y^*$ of the
problem (\ref{f5})--(\ref{f6}). Moreover, the global solution
$y_d^*$  of the reduced problem (\ref{f2.5})--(\ref{f2.7}) will
also belong to the same $\varepsilon$-neighborhood.

\begin{theorem}
\label{t1} For every problem  (\ref{f5})--(\ref{f6}) having the
unique solution $y^*$, and a given accuracy $\varepsilon>0$,
there exists a subdivision level $d$ such that in the
$\varepsilon$-neighborhood $\varepsilon(y^*)$ of the global
minimizer $y^*$ there exists a segment $h$ of the level-$d$
approximating curve $y_{d}(x)$ with the following properties:

\begin{enumerate}
\item
the segment $h$ belongs to $Q_{m+1} \cap \varepsilon(y^*)$ and has
a finite positive length;
\item
the segment $h$ contains the global minimizer $y_d^*$ of the
problem
\begin{equation}
\varphi(y_d^*)=\min\{G_{m+1}(y_d(x)):\; y_d(x)\in Q_{m+1}, x
\in[0,1]\}. \label{f3.1}
\end{equation}
\end{enumerate}
\end{theorem}

{\bf Proof.} We have assumed that the feasible region $Q_{m+1}$ is
the union of $T_{m+1}$ robust sub\-regions $Q_{m+1}^j, \; 1\le
j\le T_{m+1}$. Suppose that $y^* \in Q_{m+1}^{p}$, where
$Q_{m+1}^{p}$ is one of those subregions. Then $Q_{m+1}^{p} \cap
\varepsilon(y^*)$ is a robust set with a nonzero volume. This
means that by increasing the accuracy of approximation  it is
possible to find  an approximation level  $d$ such that the curve
$y_{d}(x)$ will have a segment $h \subset Q_{m+1}^{p} \cap
\varepsilon(y^*)$ with the required properties.

\begin{figure}[t]
\centerline{\psfig{file=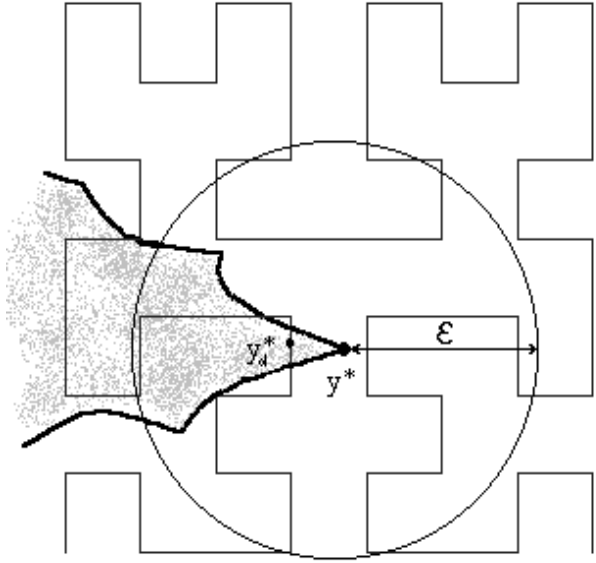,width=80mm,height=80mm,silent=yes}}
\caption{Illustration of existence of a subdivision
 level $d$ allowing an $\varepsilon$-approximation $y_{d}^*$ of
the global solution $y^*$ such that the point $y_{d}^*$ is the
global minimizer of the objective function $G_{m+1}(y_d(x))$
 on $Q_{m+1} \cap y_{d}(x)$.}
 \label{Fig_2}
\end{figure}

An illustration to Theorem~\ref{t1} is given  in Fig.~\ref{Fig_2},
where  a region $Q_{m+1}^{p}$ is shown in grey color.
Analogously, if problem  (\ref{f5})--(\ref{f6}) has more than one
solution, an approximation level $d$ providing existence of such
intervals having finite positive lengths can be found for each of
global minimizers. Of course, (\ref{f3.1}) may be satisfied only
for one interval but, anyway, all these intervals will contain
$\varepsilon$-approximations of the global solution. Thus,
Theorem~\ref{t1} allows us to concentrate the further theoretical
investigation on the behaviour of the method during the solution
of  problem (\ref{f2.10}).

\begin{lemma}
\label{l1} Let $\bar{x}$  be a limit  point of the infinite
sequence $\{x^{k}\}$ generated by the algorithm and let $i = i(k)$
be the number of an interval $[x_{i-1},x_{i}]$ containing this
point during the $k$-th iteration. Then
\begin{equation}
\lim_{k \rightarrow \infty}(x_{i(k)} - x_{i(k)-1})^{1/N} = 0,
                              \label{f3.2}
\end{equation}
and for every $\sigma  > 0$ there exists an iteration number
$n(\sigma )$, such that the  inequality
\begin{equation}
R_{i(k)} < \sigma       \label{f3.3}
\end{equation}
holds for all $k \ge  n(\sigma  )$.
\end{lemma}

{\bf Proof.} The new trial point $x^{k+1}$ from (\ref{f2.18}) or
(\ref{f2.19}) falls into an interval $[x_{i-1},x_{i}]$  (where
$i(k)=t(k)$ is determined by (\ref{f2.17})), and divides this one
into two sub-intervals
\[
[x_{i-1},x^{k+1}],\quad [x^{k+1},x_{i}].
\]
Let us show that there exists a number $\alpha$, independent of
the iteration number $k$, such that 
\begin{equation}
\max \{ x_{i}-x^{k+1}, x^{k+1}-x_{i-1} \} \le \alpha
(x_{i}-x_{i-1}), \quad 0.5\le \alpha <1,
                    \label{f3.4}
\end{equation}
holds for these intervals. In the case $\nu _{t}\neq \nu _{t-1}$
the estimate (\ref{f3.4}) is obtained immediately from
(\ref{f2.19}) by taking $\alpha_1 =0.5$. In the opposite case,
from (\ref{f2.13})--(\ref{f2.15}) it follows that
\begin{equation}
\frac{\mid z_{t}-z_{t-1}\mid ^N}{M_{t}^N} \le  x_{t}-x_{t-1},
\quad \nu _{t}=\nu _{t-1}.        \label{f3.5}
\end{equation}
From this inequality, (\ref{f2.18}), and the fact that $r > 1$ we
can conclude that in the case $\nu _{t}=\nu _{t-1}$ (\ref{f3.4})
is true for $ \alpha_2 = (r+1)/(2r)$. Since $ \alpha_2>\alpha_1$,
(\ref{f3.4}) is proved by taking $ \alpha = \alpha_2$.

Now, considering (\ref{f3.4}) together with the existence of a
sequence converging to $\bar{x}$ (this point is a limit point of
$\{x^{k}\}$), we can deduce that (\ref{f3.2}) holds. Note that in
the case when two intervals containing the point $\bar{x}$  there
exist (i.e., when $\bar{x}\in \{x^{k}\}$) the number $i=i(k)$ is
juxtaposed to the interval for which (\ref{f3.2}) takes place.

Let us prove (\ref{f3.3}). From (\ref{f2.16}) and (\ref{f3.5}) we
obtain, for intervals with indexes  $\nu _{i}=\nu _{i-1}$,
\[
R_{i} \le  (1+r^{-2})(x_{i} - x_{i-1})^{1/N} \le 2(x_{i} -
x_{i-1})^{1/N},
\]
taking into account that $\Phi (x)\ge 0$ and $r>1$. The last
estimate holds also for intervals with $\nu _{i} \ne \nu _{i-1}$.
This inequality and (\ref{f3.2}) lead to (\ref{f3.3}).

\begin{theorem}
\label{t2} Let $x^{*}$ be any solution to the problem
(\ref{f2.10})  and $j=j(k)$ be the number of an interval
$[x_{j-1},x_{j}]$ containing this point during the $k$-th
iteration. Then, if for $k\ge k^{*}$ the condition
\begin{equation}
rM_{j}>\left\{ \begin{array}{ll}
        2^{1-1/N}C_{j}+\sqrt{2^{2-2/N}C_{j}^{2}-D_{j}^{2}}, & \nu_{j-1}=\nu_{j},\\
          2C_{j},                         & \nu_{j-1}\neq \nu_{j},
               \end{array}
       \right.          \label{f3.6}
\end{equation}
 holds, the point $x^{*}$ will be a limit point of
$\{x^{k}\}$. The values $C_{j}$ and $D_{j}$ used in (\ref{f3.6})
are
\begin{equation}
C_{j}=\left\{ \begin{array}{ll}
        z_{j-1}/(x^{*}-x_{j-1})^{1/N}, & \nu_{j-1}>\nu_{j},\\
        \max \{ z_{j-1}/(x^{*}-x_{j-1})^{1/N}, z_{j}/(x_{j}-x^{*})^{1/N} \},&  \nu_{j-1}=\nu_{j},\\
        z_{j}/(x_{j}-x^{*})^{1/N}, & \nu_{j-1}<\nu_{j},
               \end{array}
       \right.          \label{f3.7}
\end{equation}
\begin{equation}
D_{j}=\left\{ \begin{array}{ll}
    \mid z_{j}-z_{j-1}\mid /(x_{j}-x_{j-1})^{1/N}, & \nu_{j-1}=\nu_{j},\\
        0,    &        \mbox{\rm otherwise}.
              \end{array}
      \right.   \label{f3.8}
\end{equation}

\end{theorem}

{\bf Proof.} Consider the case $\nu _{j-1}=\nu _{j}$. Due to
(\ref{f3.7}), we can write
\begin{equation}
z_{j-1}\le  C_{j}(x^{*}-x_{j-1})^{1/N},   \label{f3.9}
\end{equation}
\begin{equation}
z_{j}\le  C_{j}(x_{j}-x^{*})^{1/N}.           \label{f3.10}
\end{equation}
Now, by using (\ref{f3.9}), (\ref{f3.10}), and the designation
\[
\beta = (x^{*} - x_{j-1} ) / ( x_{j} - x_{j-1} ),
\]
we deduce
\[
\begin{array}{ll}
z_{j-1} + z_{j} & \leq C_{j} ( ( x^{*} - x_{j-1} )^{1/N} + ( x_{j}
- x^{*} )^{1/N} ) \\
 & =C_{j} ( \beta^{1/N} + ( 1 - \beta )^{1/N} ) ( x_{j} - x_{j-1}
)^{1/N} \\
& \leq  C_{j} ( x_{j} - x_{j-1} )^{1/N} \max \{ \beta^{1/N} + (1 -
\beta )^{1/N} : 0 \leq \beta \leq 1 \} \\
& =  2^{1-1/N} C_{j} ( x_{j} - x_{j-1} )^{1/N}.
 \end{array}
\]
By using this estimate and (\ref{f3.6}), (\ref{f3.8}), we obtain
\[
\ba{ll} R_{j(k)} & =( x_{j} - x_{j-1} )^{1/N} +
\displaystyle{\frac{( z_{j-1}\!-\! z_{j} )^{2}}{ ( r M_{j} )^{2}
(x_{j}\!-\!x_{j-1} )^{1/N}} - \frac{2( z_{j-1}\!+\!z_j)}{r M_{j}}}\\
\\
& \ge (x_j\!-\!x_{j-1})^{1/N} (1\!+\!D^{2}_{j} (r M_{j}
)^{-2}\!-\!2^{2-1/N} C_j(r M_{j} )^{-1}).
                            \ea
\]
Now, due to (\ref{f3.6}), we can conclude that
\begin{equation}
R_{j(k)} > 0.        \label{f3.11}
\end{equation}
In the case $\nu _{j-1}> \nu _{j}$ the estimate (\ref{f3.9})
takes place due to (\ref{f3.7}). From (\ref{f2.16}) it follows
that
\[
R_{j(k)} \ge  2(x_{j}-x_{j-1})^{1/N}(1-2C_{j}(rM_{j})^{-1}),
\]
and, consequently, by taking into consideration (\ref{f3.6}), the
inequality (\ref{f3.11}) holds in this case also. Truth of
(\ref{f3.11}) for $\nu _{j-1}< \nu _{j}$ is demonstrated by
analogy.

Assume now that $x^{*}$ is not a limit point of the sequence
$\{x^{k}\}$. Then, there exists a number $K_s$ such that for all
$k \ge K_s$ the interval $[x_{j-1},x_{j}]$, $j=j(k)$, is not
changed, i.e., new points will not fall into this interval.

Consider again the interval $[x_{i-1},x_{i}]$ from Lemma~\ref{l1}
containing a limit point $\bar{x}$. It follows from (\ref{f3.3})
that there exists a number $K_p$ such that
\[
R_{i(k)} < R_{j(k)}
\]
for all $k \ge  k^{*}=\max \{K_s,K_p\}$. This means that, starting
from the iteration $k^{*}$, the characteristic of the interval
$[x_{i-1},x_{i}], i = i(k), k\ge k^{*}$, is not maximal. Thus, a
trial will fall into the interval $[x_{j-1},x_{j}]$. But this fact
contradicts our assumption that $x^{*}$ is not a limit point.

\begin{corollary}
 \label{c1}
Given the conditions of Theorem~\ref{t2}, Theorem~\ref{t1} ensures
that trial points generated by the algorithm will fall in the
feasible sub-region $Q_{m+1} \cap \varepsilon(y^*)$.
\end{corollary}

Theorems~\ref{t1},~\ref{t2} generalize for the constrained case
(\ref{f2.5})--(\ref{f2.7}) results  established in~\cite{Se95}
for the Lipschitz global optimization problems with box
constraints. Usually, Lipschitz global optimization algorithms
need an overestimate (adaptive or a priori given) of the {\it
global} Lipschitz constant for the {\it whole} search region
(see, for example, survey \cite{HanJau} and references given
therein). The new algorithm does not need knowledge of the
precise Lipschitz constants for the functions $G_i(y),\, 1\le
i\le m+1,$ over the whole search region. In fact, the point $y^*$
(as all points in $D$) has up to $2^N$ images on the Peano curve.
To obtain an $\varepsilon$-approximation of $y^*$ it is
sufficient the fulfillment of the condition (\ref{f3.6}) (which is
significantly weaker than the Lipschitz condition) for only one
image of $y^*$ on a segment $h$ of the curve such that $h \subset
Q_{m+1} \cap \varepsilon(y^*)$. In the rest of the region $D$ the
constants $L_i,\, 1\le i\le m+1,$ can be underestimated.

In the preceding it has been established that, if condition
(\ref{f3.6}) is satisfied, the trial sequence $\{x^{k}\}$
generated by the algorithm converges to the global minimizer $x^*$
of the function  $\Phi (y)$ and, consequently, the corresponding
sequence $\{y_d(x^{k})\}$ converges to the global minimizer $y^*$
of the function  $\varphi (y)$. In the following theorem it is
proved that for any problem (\ref{f5})--(\ref{f6}) there exists a
continuum of values of $r$ satisfying (\ref{f3.6}).

\begin{theorem}
\label{t3} For any problem (\ref{f5})--(\ref{f6}) a value $r^{*}$
exists such that, for all $r > r^{*}$, condition (\ref{f3.6}) is
satisfied for all iteration numbers $k \ge 1$.
\end{theorem}

{\bf Proof.} Let us choose
\begin{equation}
r^{*} = 2^{3-1/N}\xi^{-1}\sqrt{N+3} \; \max \{ L_i:\; 1\le i\le
m+1\}, \label{f3.12}
\end{equation}
and take a value $r>r^{*}$. Since, due to (\ref{f2.7}),
\[
H_i = 2\, \sqrt{N+3}\; L_i,\quad 1\le i\le m+1,
\]
the value $r^{*}$ from (\ref{f3.12}) can be rewritten as
\[\begin{array}{ll}
r^{*} &= 2^{2-1/N}\xi^{-1} H,\\[5pt]
H     &= \max \{ H_i: 1\le i\le m+1\}. \end{array}
\]
For all $k \ge 1$, for the estimates $M_{j(k)}$ and the values
$C_{j(k)}$ from (\ref{f3.7}) the inequalities
\[
\xi \le M_{j(k)}, \hspace{5mm}C_{j(k)} \le H,
\]
hold. Then, when the interval $[x_{j-1},x_{j}]$ has $\nu_j =
\nu_{j-1}$, it follows from (\ref{f3.12}) that
\[
\begin{array}{lll}
rM_{j(k)} &> r^{*}M_{j(k)} &\ge r^{*}\xi = 2^{2-1/N}H \\[5pt]
          &> 2^{2-1/N}C_{j} &\ge
2^{1-1/N}C_{j} + \sqrt{2^{2-2/N}C_{j}^{2}-D_{j}^{2}}.
\end{array}
\]
To complete the proof it is sufficient to note that in the case
$\nu_j \ne \nu_{j-1}$ the last estimate can be substituted by
\[
2^{2-1/N}C_{j} > 2\,C_{j},
\]
because $N \ge 2$. Thus, for the chosen $r>r^{*}$, condition
(\ref{f3.6}) is satisfied.

Theorem~\ref{t2} establishes sufficient convergence conditions for
the method. Theorem~\ref{t3} ensures that there exists a
continuum of values of the parameter $r$ satisfying these
conditions. However, these theorems just prove {\it existence} of
such values and  can not be considered as an instrument for
determining the value $r^{*}$. This value is problem-dependent
and can not be found without  additional information about the
objective function and constraints. In presence of such
information it is possible to solve the problem (S1), otherwise
only the problem (S2) can be faced.

Condition (\ref{f3.6}) gives us a suggestion how to choose a
reasonable value of $r$ to start optimization in the case (S2).
Note that the value $D_{j}$ from (\ref{f3.9}) can be equal to
zero. The role of $M_j$ is to estimate $C_j$. Thus, if $D_{j}=0$
and $M_j=C_j$, it follows from (\ref{f3.6}) that $r$ should be
greater than $2^{2-1/N}$. Some practical advises for the choice of
the method parameters will be given in the next section.


\section{Numerical experiments}

This section presents numerical experiments   that investigate
the performance of the proposed algorithm in solving some test
problems.  The first three problems are taken from the
literature; the other problems are proposed by the authors. All
the methods have been implemented in MATLAB~\cite{Matlab} and the
experiments have been executed at a PC with Pentium III $733$MHz
processor.

{\bf Experiment 1.} In the first experiment the new algorithm is
compared with  two methods: (i) the information algorithm for
solving problems with box constraints proposed in~\cite{St85},
combined with the penalty approach; (ii) the original index
information algorithm proposed by Strongin and Markin
in~\cite{St85}--\cite{StronMark}, that does not use the local
tuning. These methods have been chosen for comparison because all
of them have similar computational cost for a single iteration (of
course, if the cost required for the search of the penalty
coefficient allowing to solve the problem is not taken in
consideration for the method (i)), they use Peano curves, and have
the same stopping rule.  Since the penalty approach needs the
objective function and constraints defined over the whole search
region, test problems enjoying this property have been chosen.

\begin{figure}[ht]
\centerline{\psfig{file=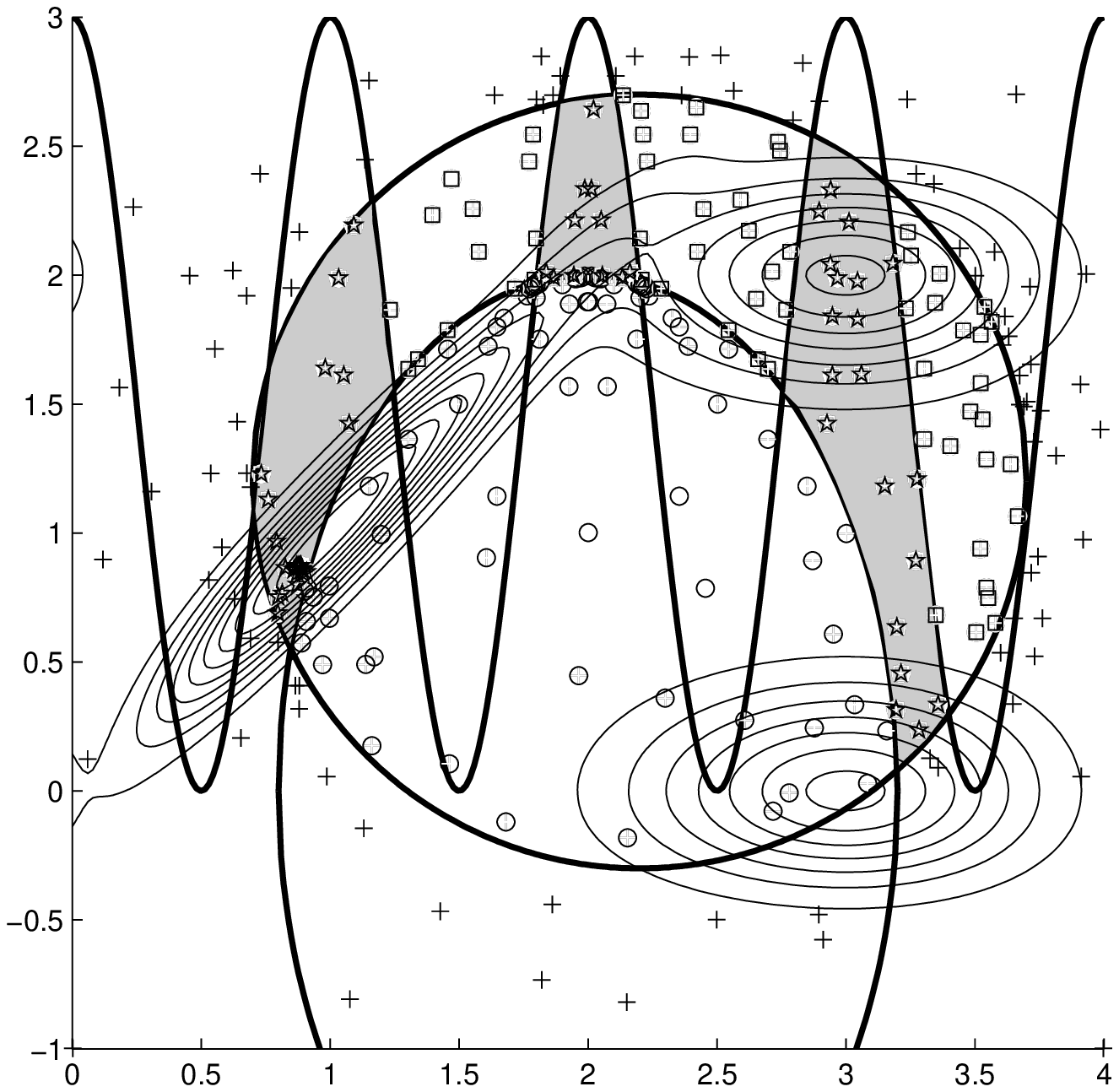,width=\textwidth,silent=yes}}
\caption{Trial points generated by  the new algorithm with $\delta
= 10^{-3}$ for Problem~1}
 \label{Fig_3}
\end{figure}

{\it Problem 1.} The first problem (see~\cite{StronMark}) is to
minimize the function
\begin{equation} \begin{array}{ll}
\varphi(y) = &
-1.5 \, y_1^2\exp(1-y_1^2-20.25\,(y_1-y_2)^2)\\[3pt]
 & - \, (0.5\,(y_1-1)(y_2-1))^4\exp(2-(0.5\,(y_1-1))^4-(y_2-1)^4)
\end{array}  \label{P1}
\end{equation}
over the rectangle $D = [0,4]\times[-1,3]$, under the constraints

\begin{displaymath} \begin{array}{l}
G_1(y) = 0.01 \,((y_1-2.2)^2 + (y_2-1.2)^2 - 2.25) \le 0,\\[3pt]
G_2(y) = 100 \,(1-(y_1-2)^2/1.44 - (0.5\,y_2)^2) \le 0,\\[3pt]
G_3(y) = 10 \,(y_2 - 1.5 - 1.5\sin(6.283\,(y_1-1.75))) \le 0.
\end{array} \end{displaymath}
The feasible region (see Fig.~\ref{Fig_3}) is the intersection of
the areas inside the circle $G_1(y) = 0$, outside the ellipses
$G_2(y) = 0$, and below the sinusoide $G_3(y) = 0$. It consists of
three disjoint and non-convex pieces with non-smooth boundaries
shown in Fig.~\ref{Fig_3} by the gray color. Here the solution is
$y^*=(0.942,0.944)$, giving the value
$\varphi(y^*)=-1.489$~\cite{SS2000}.

\begin{figure}[t]
\centerline{\psfig{file=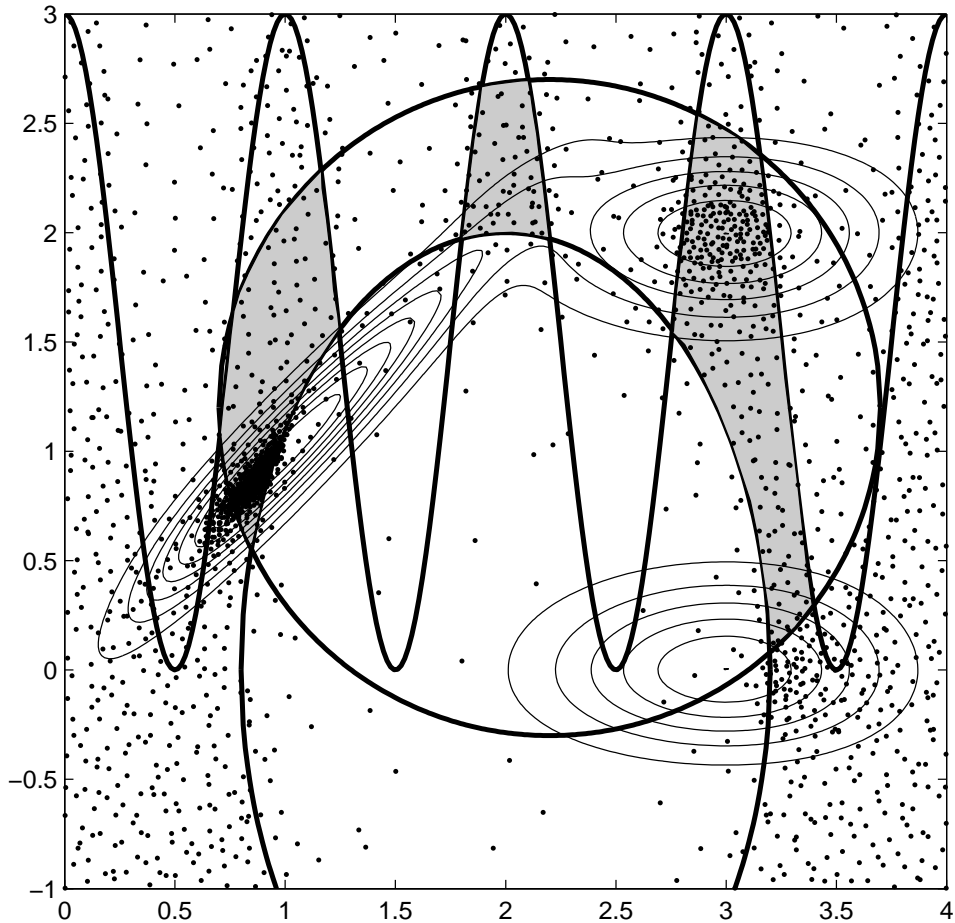,width=\textwidth,silent=yes}}
\caption{Trial points generated  by  the information algorithm
combined with the penalty approach with $\delta = 10^{-3}$ for
Problem~1}
 \label{Fig_5}
\end{figure}

\begin{figure}[t]
\centerline{\psfig{file=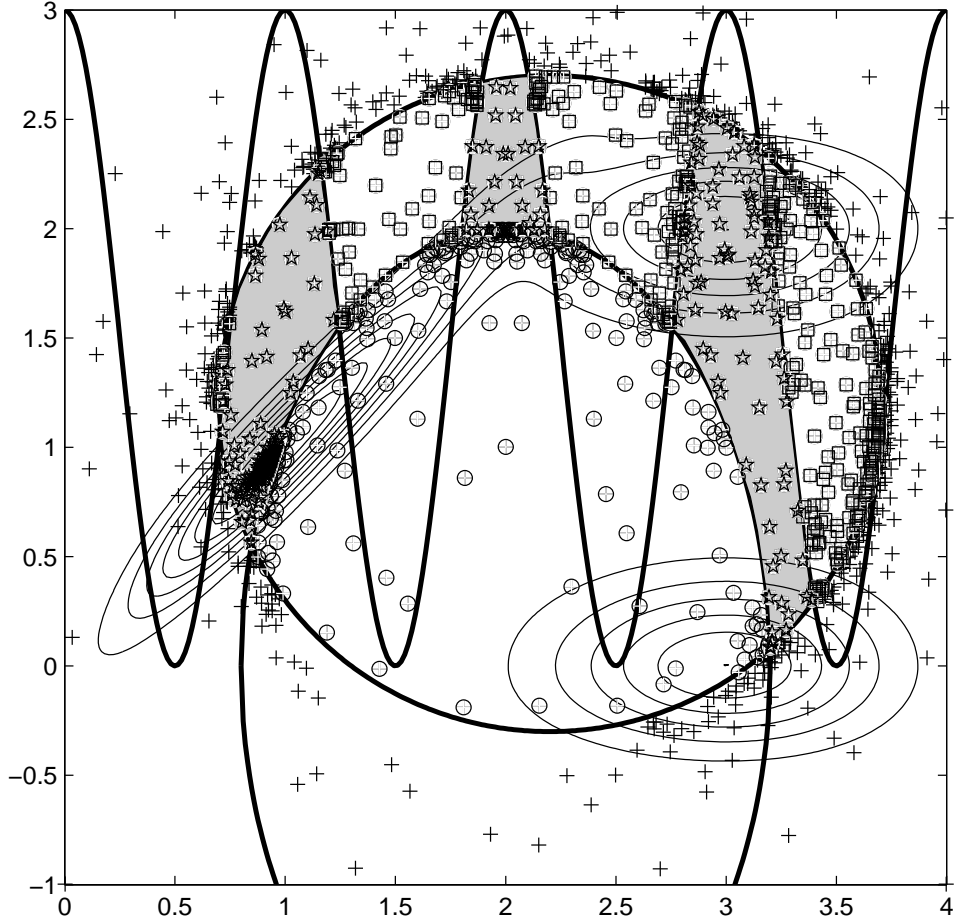,width=\textwidth,silent=yes}}
\caption{ Trial points generated by  the Strongin-Markin
algorithm with $\delta = 10^{-3}$ for Problem~1}
 \label{Fig_4}
\end{figure}

{\it Problem 2.} In this problem (see~\cite{St92})  the same
function as in the first experiment, over the same  domain $D$, is
minimized subject to the constraints
\begin{equation}  \begin{array}{l}
G_1(y) = -(y_1-2.2)^2-(y_2-1.2)^2 \le -1.21,\\[3pt]
G_2(y) = (y_1-2.2)^2+(y_2-1.2)^2 \le 1.25,
\end{array}  \label{P2}
 \end{equation}
that define a narrow annulus centered at the point $C \equiv
(2.2,1.2)$. Here the solution is $y^*=(1.088,1.088)$, giving the
value $\varphi(y^*)=-1.477$~\cite{SS2000}.

{\it Problem 3.} This is a minimization problem with four
constraints (see Section 9.2 in~\cite{Himmelblau}); the
expression of the objective function is reported in the appendix;
the search region is $D = [0,80] \times [0,80]$ and the
constraints are
\begin{equation}  \begin{array}{l}
G_1(y) = 450 - y_1 y_2 \le 0,\\[3pt]
G_2(y) = (0.1\,y_1-1)^2-y_2\le 0,\\[3pt]
G_3(y) = 8\,(y_1-40)-(y_2-30)(y_2-55)\le 0,\\[3pt]
G_4(y) = (y_1-35)(y_1-30)/125+y_2-80\le 0.
\end{array}  \label{P3}
\end{equation}
The feasible region for this problem is connected but not convex and has a
non-smooth boundary. The solution is the point
$y^*=(77.19,64.06)$, and the value of the function at this point
is $\varphi(y^*)=-59.59$~\cite{SS2000}.

{\it Problem 4.} The fourth test problem, proposed by the
authors, is defined by the search domain $D = [0,2\,\pi] \times
[0,2\,\pi]$, the non-smooth objective function
\begin{equation}
\varphi(y) = - |\sin(y_1) \sin(2\,y_2)| + 0.01\,(y_1\,y_2 +
(y_1-\pi)^2 + 3\,(y_2-\pi)^2)\,,    \label{P4}
\end{equation}
and the non-smooth constraints
\begin{displaymath}  \begin{array}{l}
G_1(y) = 1 - y_2 + \pi/2 - |\sin(2\,y_1)| + y_1/3 \le 0,\\[3pt]
G_2(y) = y_2 - 3\,\pi/2 + 4\,|(\sin(y_1+\pi)|+y_1/3-1.9 \le 0.
\end{array} \end{displaymath}
The feasible region is made of two disjoint pieces,
each of which is non-convex and has a non-smooth boundary. The
best solution, found by a fine grid over the search region, is
$y^*=(1.247,2.392)$, yielding $\varphi(y^*)= - 0.864$.

\newcommand{\xx}{\hspace*{-1mm}}

\begin{center}
\begin{table}[t]\centering
\caption{Results for Experiment 1 with $\delta = 10^{-3}$} \label{table_1}\footnotesize
\begin{tabular}{|c|rrrrrr|rrrrrr|}\hline
        & \multicolumn{6}{c|}{Strongin-Markin method} & \multicolumn{6}{c|}{New method}\\
\xx Problem \xx\xx & $n_1$ & $n_2$ & $n_3$ & $n_4$ & $n_\varphi$ & $\varphi_c$ & $n_1$ & $n_2$ & $n_3$ & $n_4$ & $n_\varphi$ & $\varphi_c$ \\[1pt] \hline
 1      &  1712 &\ \@ 953&\ \@ 673 &    -   &   279 &     -1.489  &   289 &   201 &   119 &  -     &\@\@60 &      -1.455 \\
 2      &  1631 &   672  &    -    &    -   &   246 &     -1.477  &   316 &   142 &   -   &  -     &    31 &      -1.438 \\
 3      &   351 &   348  &   316   &    140 &    88 &     -59.59  &   106 &   103 &    90 &     65 &    36 &      -59.34 \\
 4      &   771 &   359  &     -   &    -   &   163 &     -0.864  &   115 &    60 &    -  &  -     &    29 &      -0.863 \\
\hline
\end{tabular}
\end{table}
\end{center}

\begin{center}
\begin{table}[t]\centering
\caption{Results for Experiments 1 with $\delta = 10^{-4}$}\label{table_2} \footnotesize
\begin{tabular}{|c|rrrrrr|rrrrrr|}\hline
           & \multicolumn{6}{c|}{Strongin-Markin algorithm} & \multicolumn{6}{c|}{New method}\\
\xx Problem \xx\xx & $n_1$ & $n_2$ & $n_3$ & $n_4$ & $n_\varphi$ & $\varphi_c$ & $n_1$ & $n_2$ & $n_3$ & $n_4$ & $n_\varphi$ & $\varphi_c$ \\[1pt] \hline
 1      &  4494 &  2636 &  1487 &    -  &   826 &     -1.489 &   362 &   274 &   162 &    -   &   102 &      -1.489 \\
 2      &  4926 &  1941 &   -   &   -   &   781 &     -1.478 &   376 &   170 &    -   &   -    &    44 &      -1.439 \\
 3      &  1073 &  1070 &  1032 &   655 &   359 &     -59.59 &   115 &   112 &    99 &    74 &    45 &      -59.34 \\
 4      &  1867 &   787 &   -   &   -   &   355 &     -0.864 &   135 &    72 &   -    &   -    &    40 &      -0.863 \\ \hline
\end{tabular}
\end{table}
\end{center}

In all the experiments the maximum number of iterations was set to
5000, the level of approximation of the Peano curve 
to $d = 10$, the value  $\xi$ from Step 4 to $\xi = 10^{-8}$, and
the value of the reliability parameter to $r = 2.2$.  Only the
points $0$ and $1$ have been used as starting points.

In the experiments with the penalty approach the constrained
problems were reduced to the unconstrained ones as follows
\[ \varphi_{P}(y)=\varphi(y) + P \cdot
\max\left\{G_1(y),G_2(y),\dots,G_{m}(y),0\right\},  y \in D.
 \]
The coefficient $P$ has been computed by the rules:
\begin{description}
\item
-- the coefficient $P$ has been  chosen equal  to $0.1$ for all
the problems and it has  been checked whether the found solution
for each problem belongs or not to the feasible subregions;
\item
-- if it does not belong to the feasible subregions, the
 coefficient $P^*$ has been iteratively increased by $0.1$
  until a feasible solution has been found.
\end{description}

Figures \ref{Fig_5}--\ref{Fig_3} present trial points generated by
the three methods during Experiment~1 with $\delta = 10^{-3}$.
Figure~\ref{Fig_5} shows the trial points generated by the penalty
approach. Remind that both the objective function and all the
constraints have been evaluated at every point.
Figures~\ref{Fig_4} and~\ref{Fig_3} show the points where the
algorithms using the index approach have chosen to evaluate the
constraints and eventually the target function. Points where the
target function has been computed are denoted by stars; those
where only the first constraint has been evaluated (and found
violated) are denoted by crosses; circles denote points where the
first constraint was satisfied and the second was not; squares
denote points were the first two constraints were satisfied and
the third was not.

\begin{center}
\begin{table}[t]\centering
\caption{Results obtained by the information algorithm combined
with the penalty approach  Experiment 1 }\label{table_x}
\footnotesize
\begin{tabular}{|c|l|rr|rr|}\hline
Test Problem &  $P^*$  & \multicolumn{2}{c|}{$\delta = 10^{-3}$} &
\multicolumn{2}{c|}{$\delta = 10^{-4}$}\\
& & $n$ & $\varphi_c$ & $n$ & $\varphi_c$ \\
\hline
      1      & $0.1$ & $2298$       & $-1.489$ & ${\bf 5000}$ & $-1.489$ \\
      2      & $0.3$ &  $428$       & $-1.478$ & $1621$       & $-1.478$ \\
      3      & $0.8$ & ${\bf 5000}$ & $-59.59$ & ${\bf 5000}$ & $-59.59$ \\
      4      & $0.5$ &  $283$       & $-0.864$ &  $789$       & $-0.864$ \\ \hline
\end{tabular}
\end{table}
\end{center}

\begin{center}
\begin{table}[t]\centering
 \caption{Speed-up obtained in Experiment 1 for Problems 1--4}\label{table_4} \footnotesize \centering
\begin{tabular}{|c|rrrr|rrrr|}\hline
Problem       &  1  &   2  &  3  &   4  &   1  &   2  &  3  &   4    \\
    & \multicolumn{4}{c|}{$\delta = 10^{-3}$} & \multicolumn{4}{c|}{$\delta = 10^{-4}$}  \\
    \hline
  &&&&&&&&  \\
$S_1$   & 4.29 & 5.13 & 3.74 & 4.03 & 12.25 & 10.84 & 14.68 &  8.83  \\
$S_2$   & 4.65 & 7.94 & 2.44 & 5.62 &  8.10 & 17.75 &  7.98 &  8.88  \\
$S_3$   & 5.41 & 5.21 & 3.11 & 6.34 & 10.49 & 12.96 &  9.41 &
11.27  \\
$\hat{S}_1$   & 7.95 & 1.35 & 47.17 & 2.46 & 13.81 & 4.31 & 43.48 &  5.84  \\
$\hat{S}_2$   & 38.30 & 13.81 & 138.89 & 9.76 &  49.02 & 36.84 &  111.11 &  19.73  \\
$\hat{S}_3$   & 13.74 & 2.63 & 62.50 & 4.16 & 22.47 & 8.24 & 56.18
& 9.58  \\
 \hline
\end{tabular}
\end{table}
\end{center}

Table \ref{table_1} compares performance of the proposed
algorithm to the Strongin-Markin  me\-thod for $\delta = 10^{-3}$.
The data reported are the computed solution $\varphi_c$, the
number of evaluations of each constraint $n_i$, and of the target
function $n_\varphi$. To investigate the sensitivity of both
algorithms to the value of $\delta$ the tests have been repeated
with $\delta = 10^{-4}$; results are reported in
Table~\ref{table_2}. Finally, Table~\ref{table_x} presents results
obtained by  the information algorithm combined with the penalty
approach on Problems $1-4$, where $n$ is the number of iterations
executed by the algorithm. Cases where the algorithm was not able
to stop in 5000 iterations are shown in bold.

Table \ref{table_4} illustrates acceleration  reached by the new
algorithm in comparison with the two methods tested in
Experiment~1. Three speed-up indexes have been used: the first
index is the ratio $S_1$ of iterations $n_1$ of the
Strongin-Markin method over the new one; the second is the ratio
$S_2$ of $n_\varphi$ of the objective function evaluations;
finally, $S_3$ is the ratio of the sums $n_\varphi + \sum_{k=1}^m
n_k$, being the number of  summary evaluations of the objective
function and constraints. The same quantities for the penalty
approach are indicated by $\hat{S}_i, 1 \le i \le 3$. It can be
seen from Table~\ref{table_4} that the new method significantly
outperforms both methods tested. Moreover, the computational
burden increases very slowly for decreasing $\delta$ for it,
while the Strongin-Markin algorithm and the penalty approach show
a fast increase of iterations when a better accuracy is required.

Since Table \ref{table_4} shows an evident superiority of the new
method over the penalty algorithm, the latter is not used in the
further experiments (also because it is difficult to determine the
correct penalty coefficients for problems used in these
experiments).

{\bf Experiment 2.} This experiment aims at investigating the
sensitivity of both algorithms to the value of the reliability
parameter $r$. We have generated twenty pairs of random numbers
in the range $(-1,1)$, and have added any of them to the center
$C \equiv (2.2,1.2)$ of the annulus of Problem~2, thus obtaining
twenty new problems, whose solutions have been found with a fine
grid over the domain $D$.

They have been solved with both methods for either $r = 2.5$ and
$r = 3$. With $r = 2.5$ the proposed algorithm has found the
solution in 18 cases out of 20, the Strongin-Markin algorithm was
successful  in 19 cases out of 20. Note that the proposed
algorithm has found the solution in the case the Strongin-Markin
algorithm has failed. With $r = 3$, both algorithms have found
the solution in all cases. The results (number of successes,
average number of computation of the constraints and of the
objective function) are reported in Table~\ref{table_3}. With $r
= 3$, in one case the Strongin-Markin algorithm has reached the
maximum number of iterations, but anyway it has found a good
approximation to the solution. It can be seen from
Table~\ref{table_3} that the computational burden increases very
slowly for increasing $r$ for the proposed algorithm, while the
Strongin-Markin algorithm  shows a fast increase of iterations.

\begin{center}
\begin{table}[t]\centering
 \caption{Results for Experiment 2}\label{table_3} \footnotesize
\begin{tabular}{|c|cc|cc|}\hline
                  & \multicolumn{2}{c|}{$r = 2.5$} & \multicolumn{2}{c|}{$r = 3$}\\
                  & Strongin-Markin method& New method & Strongin-Markin method & New method \\[1pt] \hline
Successes         & $19/20$  & $18/20$    & $20/20$  & $20/20$  \\[3pt]
$n_1$ average        &   $2177$ &   $359$    & $2939$   &   $651$  \\
$n_2$ average        &    $874$ &   $179$    & $1150$   &   $330$  \\
$n_\varphi$ average  &    $282$ &    $46$    & $312$    &    $88$
\\
$S_1$   & -- & 6.06 & -- & 4.51 \\
$S_2$   & -- & 6.13 & -- & 3.55 \\
$S_3$   & -- &   5.71 & -- & 4.12 \\ \hline
\end{tabular}
\end{table}
  \end{center}

{\bf Experiment 3.} The last experiment involves a problem  having
changeable dimension $N$ and three constraints. The  objective
function is
\begin{equation} \varphi_{N}(y) = y_1 +
e^{\rho_{N}-|\rho_{N}^2-5\,\rho_{N}+4|}, \hspace{5mm} \rho_{N}=
(\sum_{i=1}^{N}y_i^2)^\frac{1}{2}.
 \label{P5}
  \end{equation}
The search region is $D = [-2,8] \times [-6,4]\times \ldots \times
[-6,4]$ and the constraints are
\begin{displaymath}  \begin{array}{l}
G_1(y) = (y_1-2.5)^2-6.25 + \sum_{i=2}^{N}y_i^2 \le 0,\\[3pt]
G_2(y) = -(y_1-2)^2 + 2.25 - \sum_{i=2}^{N}y_i^2 \le 0,\\[3pt]
G_3(y) = y_2 - 3\,\pi/2 + 4\,|(\sin(y_1+\pi)|+y_1/3+0.8 \le 0.\\[3pt]

\end{array}
 \end{displaymath}
Numerical experiments have been executed with  $2 \le N \le 6$
with different values of $r$ and  $\delta$. The level of
approximation of the Peano curve was set to $d = 10$ for $2 \le N
\le 5$ and to $d = 8$ for $ N =6$, the value  $\xi = 10^{-8}$. In
all the experiments the maximum number of iterations was set to
$40000$. The global solution  is $y^*=(0,\ldots,0)$, yielding
$\varphi_{N}(y^*)= \exp(-4) \simeq 0.0183$.

The obtained results are shown in Tables~\ref{table_5}
and~\ref{table_6}. In all the experiments the new algorithm is
significantly faster. It can be seen from the Table~\ref{table_6}
that starting from dimension $N=3$ the Strongin-Markin algorithm
was not able to stop in less than $40000$ iterations.

Let us comment the choice of parameters of the new method. As it
has  been mentioned in the Introduction, the statement of the
problem does not allow to create any deterministic or stochastic
algorithm that, after a finite number of evaluations of
$\varphi(y)$,  is able to provide an underestimate of the global
minimum $\varphi(y^*)$. These means  that it is impossible to say
without some additional information about the problem  which
values of parameters ensure that the method finds the global
minimizer. However, the executed experiments can advice some
recommendations.

First of all, the parameters $d$ and $\xi$ are quite technical
and can be chosen easily. To have a better accuracy of the
solution it is necessary to have a better approximation of the
Peano curve determined by $d$. Since two points in the
multidimensional space have two images at the interval $[0,1]$,
then in a computer realization (see, for example, \cite{SS2000})
for a correct work of the method these images should be
represented by two different numbers. The mapping is such that
for a dimension $N$ and a level $d$ two numbers $a,b \in [0,1]$
will be considered different if $|a-b|>2^{-dN}$. For example, in
double precision the minimal representable number is $2^{-52}$.
Thus, if one works with a simple realization of the Peano curve
using double precision, the product $dN$ should be less than
$52$. Of course, more sophisticated implementations of the Peano
curve realizing numbers with more digits allow to increase this
number.

The second parameter is $\xi$ and it is chosen equal to a small
number (see experiments). It has been introduced in the method to
ensure that it works correctly in the case  $l_i=c_i=r_i=0$ where
$l_i, c_i, r_i$ are from (\ref{f2.14}). The case when the real
H\"{o}lder constants are less than $\xi$ is degenerous for the
method because the local tuning is not used and the method works
slower.

The parameters $r$ and $\delta$ are chosen by the user on the
basis of additional information about the problem or just by
increasing $r$ and decreasing $\delta$ (see Table~\ref{table_5})
trying to obtain  a satisfactory result (see statement (S2) in
Introduction). It can be also seen from  Table~\ref{table_5} that
by increasing dimension of the problem it is necessary to change
$r$ and $\delta$ in the same way. This happens because, in spite
of increasing dimension, the reduced problem is always determined
over the interval $[0,1]$ independently on the dimension $N$. As a
consequence, the reduced problem has more  local minima, and  to
find the global one it is necessary to increase reliability of
the method (by increasing $r$) and to make the search more
accurate (by decreasing $\delta$).

\begin{table}
\centering \caption{Results for the new method on
Problems~(\ref{P5}) with dimension two to six} \label{table_5}
\footnotesize
\hspace*{-2cm}\begin{tabular}{|lr|llllll|rrrr|c|}\hline
 $r$  &   $\delta$    &    $y_1$    &   $y_2$    &    $y_3$ & $y_4$    &    $y_5$  &   $y_6$   &    $n_1$   &   $n_2$ &
$n_3$    &    $n_{\varphi}$     &     $\varphi_c$ \\
 \hline
2.35 &  $10^{-3}$  &   0.0034 & 0.0596 &    -    &     - & - & - &
130   &   51  &   41  &   36 &
   0.0295 \\
2.35& $10^{-4}$  & 0.0007 &  0.0596 &    -   &      -     & - &
-   &       147 &     57   &  47   &  42  & 0.0268\\
2.45 & $10^{-3}$ &   0.0068 &  0.1185& 0.1182 &    -    & - & - &
819   &
470 &   377&   249  & 0.0555\\
2.45 & $10^{-4}$  & 0.0052& -0.0088 & 0.0107 &  -     & - & - &
5847 &  4723 &  3975 &
3313 &  0.0252\\
2.7  & $10^{-3}$&      0.2540& -0.1553& -0.0674 & 0.4600  & - & -
&      708    &  533 &   364 & 249 &
0.6247\\
2.7  & $10^{-4}$  & 0.0264 & -0.0771& 0.0400 & -0.0679 & - &
-   &  14441& 12485&10081 & 7410 & 0.0621\\
3.3  & $10^{-3}$ & 0.5146& -0.1650& -0.0771& -0.0924& -0.3799 & -
&   620 & 206   &164 &
111& 1.1704\\
3.3  & $10^{-4}$   &    0.0654 & 0.0303& -0.2920& 0.0596&
-0.0692  &    - & 18535 & 12289 & 10761 & 7344 & 0.1748\\
3.35 &$10^{-3}$ &   0.5977 & -0.1992 & -0.0039& -0.3298&1.0898&
-0.6289  &
445  &    39    & 35  &  25  & 1.9579\\
3.35 &$10^{-4}$  & 0.3193 & 0.1523& -0.3945& 0.1523& 0.3477
&0.1523&  6610  &   2486 & 2292 &
1534  &  0.9688\\
3.35 & $5 \cdot 10^{-5} $ & 0.0508& -0.0430&  0.0742& -0.0430&
-0.1992& 0.0504& 19105 &  9871 &  8989 & 6233 & 0.1208 \\
\hline
\end{tabular}
\end{table}

\begin{table}
\centering \caption{Results for the Strongin-Markin algorithm on
Problems~(\ref{P5}) with dimension two to six}\label{table_6}
\footnotesize
\hspace*{-2cm}\begin{tabular}{|lr|llllll|rrrr|c|}\hline $r$  &
$\delta$    &    $y_1$    &   $y_2$    &    $y_3$ & $y_4$ & $y_5$
&   $y_6$   &    $n_1$   &   $n_2$ &
$n_3$    &    $n_{\varphi}$     &     $\varphi_c$ \\
 \hline
2.35 & $10^{-3}$ & 0.0025 &  0.0010 & - & - & - & - &  5071 &  3811 &  3314 &  3066 & 0.0211\\
2.35 & $10^{-4}$ & 0.0000 &  0.0010 & - & - & - & - & 39883 & 33277 & 31108 & 30423 & 0.0185\\
2.45 & $10^{-3}$ & 0.0388 &  0.0010 &  0.0498 &  - & - & - & 11003 &  5734 &  4911 &  3736 & 0.0654\\
2.45 & $10^{-4}$ & 0.0068 &  0.0400 & -0.0351 &  - & - & - & {\bf 40000} & 22591 & 19901 & 16023 & 0.0320\\
2.7  & $10^{-3}$ & 0.1729 &  0.0693 & -0.1162 & -0.1860 & - & - & 5067 & 1833 & 1680 & 1037 & 0.2676\\
2.7  & $10^{-4}$ & 0.0459 &  0.1210 & -0.2139 &  0.0693 & - & - & {\bf 40000} & 17699 & 15402 & 10566 & 0.1271\\
3.3  & $10^{-3}$ & 0.5635 &  0.3037 & -0.0381 &  0.3916 &  0.0034 & - & 2661 & 432 & 381 & 220 & 1.5087\\
3.3  & $10^{-4}$ & 0.1436 &  0.3330 & -0.2627 & -0.0869 &  0.0492 & - & {\bf 40000} & 10487 & 9661 & 6137 & 0.3763\\
3.35 & $10^{-3}$ & 0.4756 & -0.2920 & -0.1846 & -0.0186 & -0.4482 & -1.1709 & 650 & 23 & 22 & 18 & 1.9356\\
3.35 & $10^{-4}$ & 0.3877 &  0.1377 &  0.4111 &  0.2061 & -0.3115 & -0.1260 & {\bf 40000} & 5413 & 5091 & 3251 & 1.1450\\
3.35 & $5\cdot 10^{-5}$ & 0.3877 &  0.1377 &  0.4111 &  0.2061 & -0.3115 & -0.1260 & {\bf 40000} & 5413 & 5091 & 3251 & 1.1450\\
 \hline
\end{tabular}
\end{table}


\section{Conclusions}

In this paper,   a novel global optimization algorithm for solving
multi-dimensional Lipschitz global optimization problems with
multiextremal partially defined constraints and objective function
has been presented. The new algorithm uses Peano type
space-filling curves and the index scheme to reduce the original
constrained problem to a H\"{o}lder one-dimensional problem. Local
tuning on the behaviour of the objective function and constraints
is executed during the work of the global optimization procedure
in order to accelerate the search. The new method works without
introducing any penalty coefficients and/or auxiliary variables.
Convergence conditions of a new type  have been established for
the algorithm.

The new algorithm enjoys the following properties:
\begin{description}
\item
\hspace{8mm}{--} in order to guarantee  convergence to a global
minimizer $y^{*}$ it is not necessary to know the precise
Lipschitz constants for the objective function and constraints
for the whole search region; on the contrary, only fulfillment of
condition (\ref{f3.6}) for a segment of the Peano curve
containing one of the images  of $y^{*}$ is required;
\item
\hspace{8mm}{--} the usual problem of determining the moment to
stop the global procedure in order to start a local search does
not arise because local information is taken into consideration
throughout the duration of the global search;
\item
\hspace{8mm}{--} local information is taken  into account not only
in the neighborhood of a global minimizer but also in the whole
search region, thus permitting a significant acceleration of the
search;
\item
\hspace{8mm}{--} thanks to usage of the index scheme, the new
algorithm does not introduce any additional parameters and/or
variables. Constraints are evaluated at every point one at a time
until the first violation of one of them, after that the rest of
constraints and the objective function are not evaluated at this
point.  In its turn, the objective function is evaluated only for
that points where all the constraints have been satisfied.
\end{description}

The algorithm has been tested on a number of problems taken from
literature. Numerical results show the good performance of the new
technique in comparison with the two methods taken from
literature.

\section*{Appendix} \
The target function for Problem 3 is
\begin{displaymath}
\begin{array}{ll}
\varphi(y) = & - ( B_1 + B_2\, y_1 + B_3\, y_1^2 + B_4\, y_1^3 +
B_5\,
y_1^4 + B_6\, y_2 + B_7\, y_1 y_2 + B_8\, y_1^2 y_2\\[3pt]
& +\, B_9\, y_1^3 y_2 + B_{10}\, y_1^4 y_2 + B_{11}\, y_2^2 +
B_{12}\, y_2^3
+ B_{13}\, y_2^4 + B_{14}/(1+y_2) \\[3pt]
&+\, B_{15}\, y_1^2 y_2^2 + B_{16}\, y_1^3 y_2^2 + B_{17}\, y_1^3
y_2^3 +
B_{18}\, y_1 y_2^2 + B_{19}\, y_1 y_2^3 \\[3pt]
& +\, B_{20} \exp(0.0005\, y_1 y_2))
\end{array}
\end{displaymath}
where coefficients $B_1, \ldots, B_{20}$ are

\small \noindent $B_1=   75.1963666677;
    B_2=   -3.8112755343;
    B_3=    0.1269366345;
    B_4=   -0.0020567665;\\
    B_5=    0.0000103450;
    B_6=   -6.8306567613;
    B_7=    0.0302344793;
    B_8=   -0.0012813448;\\
    B_9=    0.0000352559;
    B_{10}=  -0.0000002266;
    B_{11}=   0.2564581253;
    B_{12}=  -0.0034604030;\\
    B_{13}=   0.0000135139;
    B_{14}= -28.1064434908;
    B_{15}=  -0.0000052375;
    B_{16}=  -0.0000000063;\\
    B_{17}0 =   0.0000000007;
    B_{18}=   0.0003405462;
    B_{19}=  -0.0000016638;
    B_{20}=  -2.8673112392.$
\normalsize

\bibliography{constr}

\bibliographystyle{plain}

\end{document}